\documentclass[useAMS]{gITR2e}
\newcommand{\bfx}{{\bf x}}
\newcommand{\bfxp}{{{\bf x}^\prime}}
\newcommand{\wbfx}{\widehat{\bf x}}
\newcommand{\wbfxp}{{\widehat{\bf x}^\prime}}
\newcommand{\N}{{\mathbf N}}

\newcommand{\Q}{{\mathbf Q}}
\newcommand{\R}{{\mathbf R}}
\newcommand{\Si}{{\mathbf S}}

\newcommand{\g}{{\mathfrak g}}
\newcommand{\li}{{\mathfrak l}}

\newcommand{\Z}{{\mathbf Z}}
\newcommand{\C}{{\mathbf C}}

\newcommand{\mcg}{{\mathcal G}}

\makeatletter
\def\eqnarray{\stepcounter{equation}\let\@currentlabel=\theequation
\global\@eqnswtrue
\tabskip\@centering\let\\=\@eqncr
$$\halign to \displaywidth\bgroup\hfil\global\@eqcnt\z@
  $\displaystyle\tabskip\z@{##}$&\global\@eqcnt\@ne
  \hfil$\displaystyle{{}##{}}$\hfil
  &\global\@eqcnt\tw@ $\displaystyle{##}$\hfil
  \tabskip\@centering&\llap{##}\tabskip\z@\cr}

\def\endeqnarray{\@@eqncr\egroup
      \global\advance\c@equation\m@ne$$\global\@ignoretrue}

\def\@yeqncr{\@ifnextchar [{\@xeqncr}{\@xeqncr[5pt]}}
\makeatother

\begin{document}
\doi{10.1080/10652460YYxxxxxxx}
 \issn{1476-8291}
\issnp{1065-2469}
 \jvol{00} \jnum{00} \jyear{2009} \jmonth{January}

\markboth{Howard S.~Cohl}{Integral Transforms and Special Functions}


\title{
On a generalization of the generating function for Gegenbauer polynomials
}

\author{Howard S.~Cohl$^{\rm a}$$^{\ast}$\thanks{$^\ast$Corresponding author. Email: howard.cohl@nist.gov
\vspace{6pt}}
\\\vspace{6pt}  $^{\rm a}${\em{Applied and Computational Mathematics Division,
National Institute of Standards and Technology,
Gaithersburg, Maryland, U.S.A.
}} 
\\\vspace{6pt}\received{v3.5 released October 2008} }

\maketitle

\begin{abstract}
A generalization of the generating function for Gegenbauer polynomials is introduced
whose coefficients are given in terms of associated Legendre functions of the second kind.
We discuss how our expansion represents a generalization of several previously derived 
formulae such as Heine's formula and Heine's reciprocal square-root identity.
We also show how this expansion can be used to compute hyperspherical harmonic
expansions for power-law fundamental solutions of the polyharmonic equation. 
\bigskip
\begin{classcode}
35A08; 35J05; 32Q45; 31C12; 33C05; 42A16
\end{classcode}\bigskip
\end{abstract}

\section{Introduction}
\label{Introduction}

Gegenbauer polynomials $C_n^\nu(x)$ are given as the coefficients of 
$\rho^n$ for the generating function $(1+\rho^2-2\rho x)^{-\nu}$.  The study of these
polynomials was pioneered in a series of papers by Leopold 
Gegenbauer (Gegenbauer (1874,1877,1884,1888,1893)
\cite{Gegenbauer1874,Gegenbauer1877,Gegenbauer1884,Gegenbauer1888,Gegenbauer1893}).
The main result which this paper relies upon is 
{\sc Theorem} \ref{geneneralizationofgeneratingfuncitonforgegenbauerpoly} below.
This theorem gives a generalized expansion over Gegenbauer polynomials 
$C_n^\mu(x)$ of the algebraic function $z\mapsto (z-x)^{-\nu}$.
Our proof is combinatoric in nature and 
has great potential for proving new expansion formulae
which generalize generating functions. Our methodology can in principle
be applied to any generating function 
for hypergeometric orthogonal polynomials, of which 
there are many (see for instance 
Srivastava \& Manocha (1984) \cite{SriManocha}; 
Erd{\'e}lyi {\it et al.}~(1981) \cite{ErdelyiHTFII}).
The concept of the proof is to start with a 
generating function and use a connection 
formula to express the orthogonal polynomial as a finite 
series in polynomials of the same type with different parameters.  
The resulting formulae will then produce new expansions for the polynomials
which result from a limiting process, e.g.,
Legendre polynomials and Chebyshev polynomials of the first and second kind.
Connection formulae for classical orthogonal polynomials
and their $q$-extensions are well-known 
(see Ismail (2005) \cite{Ismail}).
In this paper we applied this method of proof to the generating function
for Gegenbauer polynomials.  

This paper is organized as follows.
In Section \ref{GeneralizationofthegeneratingfunctionforGegenbauerpolynomials}
we derive a complex generalization of the generating function for Gegenbauer polynomials.
In Section \ref{GeneralizationsExtensionsandApplications} we discuss how our complex
generalization reduces to 
previously derived expressions and leads to extensions in appropriate 
limits.
In Section \ref{dge3} we use our complex expansion to generalize a formula originally 
developed by Sack (1964) \cite{Sacka} on $\R^3$, to compute an expansion in terms 
of Gegenbauer polynomials for complex powers of 
the distance between two points on a $d$-dimensional Euclidean space for $d\ge 2$.
\medskip

Throughout this paper we rely on the following definitions.  
For $a_1,a_2,a_3,\ldots\in\C$, if $i,j\in\Z$ and $j<i$ then
$\sum_{n=i}^{j}a_n=0$ and $\prod_{n=i}^ja_n=1$,
where $\C$ represents the complex numbers.  
The set of natural numbers is given by $\N:=\{1,2,3,\ldots\}$, the set
$\N_0:=\{0,1,2,\ldots\}=\N\cup\{0\}$, and the set
$\Z:=\{0,\pm 1,\pm 2,\ldots\}.$
The sets $\Q$ and $\R$ represent the rational and real numbers respectively.

\section{Generalization of the generating function for Gegenbauer polynomials}
\label{GeneralizationofthegeneratingfunctionforGegenbauerpolynomials}

We present the following generalization of the generating function for Gegenbauer
polynomials whose coefficients are given in terms of associated Legendre functions
of the second kind.

\begin{theorem}
Let $\nu\in\C\setminus-\N_0,$ $\mu\in(-1/2,\infty)\setminus\{0\}$
and $z\in\C\setminus(-\infty,1]$ on any ellipse with foci at $\pm 1$
with $x$ in the interior of that ellipse. Then

\begin{equation}
\frac{1}
{(z-x)^\nu}
=\frac{2^{\mu+1/2}\Gamma(\mu)e^{i\pi(\mu-\nu+1/2)}}{\sqrt{\pi}\,\Gamma(\nu)
(z^2-1)^{(\nu-\mu)/2-1/4}}
\sum_{{n}=0}^\infty({n}+\mu)
Q_{{n}+\mu-1/2}^{\nu-\mu-1/2}(z)
C_{n}^{\mu}(x).
\label{biggeneralizationgegen2}
\end{equation}
\label{geneneralizationofgeneratingfuncitonforgegenbauerpoly}
\end{theorem}
\noindent If one substitutes $z=(1+\rho^2)/(2\rho)$ in (\ref{biggeneralizationgegen2}) 
with $0<|\rho|<1$, then one obtains an alternate expression with $x\in[-1,1],$
\begin{eqnarray}
&&\hspace{-0.6cm}
\frac{1}
{(1+\rho^2-2\rho x)^\nu}
=
\frac{\Gamma(\mu)e^{i\pi(\mu-\nu+1/2)}}
{\sqrt{\pi}\,\Gamma(\nu)
\rho^{\mu+1/2}(1-\rho^2)^{\nu-\mu-1/2}
}\nonumber\\[0.2cm]
&&\hspace{3.8cm}\times\sum_{{n}=0}^\infty
({n}+\mu)
Q_{{n}+\mu-1/2}^{\nu-\mu-1/2}
\left(\frac{1+\rho^2}{2\rho}\right)
C_{n}^\mu(x).
\label{generalizationotthergegen}
\end{eqnarray}
One can see that by replacing $\nu=\mu$ in
(\ref{generalizationotthergegen}),
and using (8.6.11) in Abramowitz \& Stegun (1972) \cite{Abra}, that these
formulae are generalizations of the generating function for Gegenbauer
polynomials 
(first occurence in Gegenbauer (1874) \cite{Gegenbauer1874})
\begin{equation}
\frac{1}
{\left(1+\rho^2-2\rho x\right)^\nu}=\sum_{{n}=0}^\infty C_{n}^\nu(x) \rho^{n},
\label{gengegen}
\end{equation}
\noindent where $\rho\in\C$ with $|\rho|<1$ and 
$\nu\in(-1/2,\infty)\setminus\{0\}$ 
(see for instance (18.12.4) in Olver {\it et al.}
(2010) \cite{NIST}).
The Gegenbauer polynomials $C_n^\nu:\C\to\C$ can be defined by
\begin{equation}
C_n^\nu(x):=\frac{(2\nu)_n}{n!}\,{}_2F_1\left(-n,n+2\nu;\nu+\frac12;\frac{1-x}{2}\right),
\label{gegpolydefn}
\end{equation}
where $n\in\N_0$, $\nu\in(-1/2,\infty)\setminus\{0\},$ and
${}_2F_1:\C^2\times(\C\setminus-\N_0)\times
\{z\in\C:|z|<1\}\to\C$, the Gauss hypergeometric function, can be defined 
in terms of the following infinite series
\begin{equation}
{}_{2}F_1(a,b;c;z):=
\sum_{n=0}^\infty \frac{(a)_n(b)_n}{(c)_n}
\frac{z^n}{n!}
\label{gauss2F1}
\end{equation}
(see (2.1.5) in Andrews, Askey \& Roy 1999), and elsewhere by 
analytic continuation.
The Pochhammer symbol (rising factorial) $(\cdot)_{n}:\C\to\C$ is defined by
\[
(z)_n:=\prod_{i=1}^n(z+i-1),
\]
where $n\in\N_0$.
For the Gegenbauer polynomials $C_n^\nu(x)$, we refer to ${n}$ and $\nu$ as the degree and 
order respectively.

\medskip

\noindent {\it Proof }\ 
Consider the generating function for Gegenbauer polynomials
(\ref{gengegen}).
The connection relation which expresses a Gegenbauer polynomial with order $\nu$ 
as a sum over Gegenbauer polynomials with order $\mu$
is given by
\begin{eqnarray}
\hspace{-3.0cm}C_n^\nu(x)&=&\frac{(2\nu)_n}{(\nu+\tfrac12)_n}\sum_{k=0}^n
\frac{(\nu+k+\tfrac12)_{n-k}\,(2\nu+n)_k\,(\mu+\tfrac12)_k\,\Gamma(2\mu+k)}
{(n-k)!\,(2\mu)_k\,\Gamma(2\mu+2k)}\nonumber\\[0.05cm]
&&{}\hspace{2.5cm}\times{}_3F_2\left(
\begin{array}{c}
-n+k,n+k+2\nu,\mu+k+\tfrac12\\ 
\nu+k+\tfrac12,2\mu+2k+1
\end{array}
;1
\right)C_k^\mu(x).
\label{connection}
\end{eqnarray}
This connection relation can be derived by starting 
with {\sc Theorem}
9.1.1 in Ismail (2005) \cite{Ismail} combined with 
(see for instance (18.7.1) in Olver {\it et al.} (2010) \cite{NIST})
\begin{equation}
C_n^\nu(x)=\frac{(2\nu)_n}{(\nu+\tfrac12)_n}P_n^{(\nu-1/2,\nu-1/2)}(x),
\label{symmetricJacobi}
\end{equation}
i.e., the Gegenbauer polynomials are given as
symmetric Jacobi polynomials.  
The Jacobi polynomials
$P_n^{(\alpha,\beta)}:\C\to\C$ can be defined as
(see for instance (18.5.7) in Olver {\it et al.} (2010) \cite{NIST})
\begin{equation}
P_n^{(\alpha,\beta)}(x):=\frac{(\alpha+1)_n}{n!}\,
{}_2F_1\left(-n,n+\alpha+\beta+1;\alpha+1;\frac{1-x}{2}\right),
\label{Jacobidefn}
\end{equation}
where $n\in\N_0,$ and $\alpha,\beta>-1$
(see Table 18.3.1 in Olver {\it et al.} (2010) \cite{NIST}).
The generalized hypergeometric function
${}_3F_2:\C^3\times(\C\setminus-\N_0)^2\times
\{z\in\C:|z|<1\}\to\C$ can be defined in terms of the following infinite series
\[
{}_3F_2\left(
\begin{array}{c}
a_1,a_2,a_3\\ 
b_1,b_2
\end{array}
;z
\right):=\sum_{n=0}^\infty
\frac{(a_1)_n(a_2)_n(a_3)_n}{(b_1)_n(b_2)_n}\frac{z^n}{n!}.
\]
If we replace the Gegenbauer polynomial in the generating 
function (\ref{gengegen}) using the 
connection relation (\ref{connection}), we obtain a double summation 
expression over $k$ and $n$.  By reversing
the order of the summations (justification by Tannery's theorem) and shifting the $n$-index by $k,$ we obtain
after making some reductions and simplifications, the following 
double-summation representation 
\begin{eqnarray}
\hspace{-0.7cm}\frac{1}{(1+\rho^2-2\rho x)^\nu}&=&\frac{\sqrt{\pi}\,\Gamma(\mu)}{2^{2\nu-1}\Gamma(\nu)}
 \sum_{k=0}^\infty C_k^\mu(x)\frac{\rho^k}{2^{2k}}\frac{\mu+k}
{\Gamma(\nu+k+\tfrac12)\Gamma(\mu+k+1)}\nonumber\\[0.2cm]
&&\hspace{-0.7cm}
\times\sum_{n=0}^\infty\frac{\Gamma(2\nu+2k)}{n!}
\,{}_3F_2\left(
\begin{array}{c}
-n,n+2k+2\nu,\mu+k+\tfrac12\\ 
\nu+k+\tfrac12,2\mu+2k+1
\end{array}
;1
\right).
\label{intermediate}
\end{eqnarray}
The ${}_3F_2$ generalized hypergeometric function appearing the above formula
may be simplified using Watson's sum 
\begin{eqnarray*}
&&\hspace{-0.3cm}{}_3F_2\left(
\begin{array}{c}
a,b,c\\ 
\frac12(a+b+1),2c
\end{array}
;1
\right)\\[0.2cm]
&&\hspace{1.3cm}=\frac
{\sqrt{\pi}\,\Gamma\left(c+\frac12\right)\Gamma\left(\frac12(a+b+1)\right)
\Gamma\left(c+\frac12(1-a-b)\right)}
{\Gamma\left(\frac12(a+1)\right)\Gamma\left(\frac12(b+1)\right)
\Gamma\left(c+\frac12(1-a)\right)\Gamma\left(c+\frac12(1-b)\right)},
\end{eqnarray*}
where $\mbox{Re}(2c-a-b)>-1$
(see for instance (16.4.6) in Olver
{\it et al.} \cite{NIST}), therefore
\begin{eqnarray}
\hspace{-1.1cm}\frac{1}{\Gamma(\nu+k+\tfrac12)\Gamma(\mu+k+1)}
\,{}_3F_2\left(
\begin{array}{c}
-n,n+2k+2\nu,\mu+k+\tfrac12\\ 
\nu+k+\tfrac12,2\mu+2k+1
\end{array}
;1
\right)&&\nonumber\\[0.2cm]
&&\hspace{-8.8cm}
=\frac{\sqrt{\pi}\,\Gamma(\mu-\nu+1)}
{\Gamma\left(\tfrac{1-n}{2}\right)\Gamma\left(\nu+k+\tfrac{n+1}{2}\right)
\Gamma\left(\mu+k+1+\tfrac{n}{2}\right)\Gamma\left(\mu-\nu+1-\tfrac{n}{2}\right)},
\label{3F2inourcase}
\end{eqnarray}
for $\mbox{Re}(\mu-\nu)>-1$.
By inserting (\ref{3F2inourcase}) in 
(\ref{intermediate}), it follows that
\begin{eqnarray}
\hspace{-0.3cm}\frac{1}{(1+\rho^2-2\rho x)^\nu}&=&
\frac{\pi\Gamma(\mu)\Gamma(\mu-\nu+1)}{2^{2\nu-1}\Gamma(\nu)}
\sum_{k=0}^\infty (\mu+k)C_k^\mu(x)\frac{\rho^k}{2^{2k}}\nonumber\\[0.2cm]
&&\hspace{-1.6cm}\times\sum_{n=0}^\infty
\frac{\rho^n\Gamma(2\nu+2k+n)}
{n!\Gamma\left(\tfrac{1-n}{2}\right)
\Gamma\left(\nu+k+\tfrac{n+1}{2}\right)
\Gamma\left(\mu+k+1+\tfrac{n}{2}\right)
\Gamma\left(\mu-\nu+1-\tfrac{n}{2}\right)}.\nonumber
\end{eqnarray}
It is straightforward to show using (\ref{gauss2F1}) and
\[
\Gamma(z-n)=(-1)^n\frac{\Gamma(z)}{(-z+1)_n},
\]
for $n\in\N_0$ and $z\in\C\setminus-\N_0$,
and the duplication
formula (i.e.,~(5.5.5) in Olver {\it et al.} (2010) \cite{NIST})
\[
\Gamma(2z)=\frac{2^{2z-1}}{\sqrt{\pi}}\,\Gamma(z)\Gamma\left(z+\frac12\right),
\]
provided $2z\not\in-\N_0$,
that 
\begin{eqnarray}
\hspace{-0.2cm}\sum_{n=0}^\infty
\frac{\rho^n\Gamma(2\nu+2k+n)}
{n!\Gamma\left(\tfrac{1-n}{2}\right)
\Gamma\left(\nu+k+\tfrac{n+1}{2}\right)
\Gamma\left(\mu+k+1+\tfrac{n}{2}\right)
\Gamma\left(\mu-\nu+1-\tfrac{n}{2}\right)}&&\nonumber\\[0.2cm]
&&\hspace{-8.0cm}=\frac{2^{2\nu+2k-1}\Gamma(\nu+k)}{\pi\Gamma(\mu+k+1)\Gamma(\mu-\nu+1)}
\,{}_2F_1\left(\nu+k,\nu-\mu;\mu+k+1;\rho^2\right),\nonumber
\end{eqnarray}
so therefore
\begin{eqnarray}
\hspace{-0.8cm}\frac{1}{(1+\rho^2-2\rho x)^\nu}&=&
\frac{\Gamma(\mu)}{\Gamma(\nu)}
\sum_{k=0}^\infty (\mu+k)C_k^\mu(x)\rho^k \nonumber\\[0.2cm]
&&\hspace{0.2cm}\times\frac{\Gamma(\nu+k)}{\Gamma(\mu+k+1)}
\,{}_2F_1\left(\nu+k,\nu-\mu;\mu+k+1;\rho^2\right).
\label{almostlast}
\end{eqnarray}
Finally utilizing the quadratic transformation of the hypergeometric 
function 
\[
{}_2F_1(a,b;a-b+1;z)=(1+z)^{-a}
{}_2F_1\left(\frac{a}{2},\frac{a+1}{2};a-b+1;\frac{4z}{(z+1)^2}\right),
\]
for $|z|<1$ (see (3.1.9) in Andrews, Askey \& Roy (1999) \cite{AAR}), 
combined with the definition of the associated Legendre function of the second 
kind $Q_\nu^\mu:\C\setminus(-\infty,1]\to\C$ in terms of the Gauss 
hypergeometric function
\begin{eqnarray}
&&\hspace{-0.6cm}Q_\nu^\mu(z):=\frac{\sqrt{\pi}\,e^{i\pi\mu}\Gamma(\nu+\mu+1)(z^2-1)^{\mu/2}}
{2^{\nu+1}\Gamma(\nu+\frac32)z^{\nu+\mu+1}}\nonumber\\[0.2cm]
&&\hspace{3.8cm}\times{}_2F_1\left(
\frac{\nu+\mu+2}{2},
\frac{\nu+\mu+1}{2};
\nu+\frac32; \frac{1}{z^2}
\right),
\label{assoclegQseckinddefn2F1}
\end{eqnarray}
for $|z|>1$ and $\nu+\mu+1\notin-\N_0$ (cf.~Section 14.21 and (14.3.7) in 
Olver {\it et al.} (2010) \cite{NIST}), one can show that
\begin{eqnarray*}
&&{}_2F_1\left(\nu+k,\nu-\mu;\mu+k+1;\rho^2\right)\\[0.2cm]
&&\hspace{2.0cm}=\frac{\Gamma\left(\mu+k+1\right)
e^{i\pi(\mu-\nu+1/2)}}
{\sqrt{\pi}\,\Gamma(\nu+k)\rho^{\mu+k+1/2}(1-\rho^2)^{\nu-\mu-1/2}}
Q_{k+\mu-1/2}^{\nu-\mu-1/2}\left(\frac{1+\rho^2}{2\rho}\right),\nonumber
\end{eqnarray*}
which when used in (\ref{almostlast}) produces 
(\ref{generalizationotthergegen}).
Since the Gegenbauer polynomial is just a symmetric Jacobi polynomial
(\ref{symmetricJacobi}), through {\sc Theorem 9.1.1} in Szeg\H{o} (1959) \cite{Szego}
(Expansion of an analytic function in a Jacobi series), since $f_z:\C\to\C$ 
defined by $f_z(x):=(z-x)^{-\nu}$
is analytic in $[-1,1]$, then the above expansion in Gegenbauer 
polynomials is convergent 
if the point $z\in\C$ lies on any ellipse with foci at $\pm 1$ and $x$ can lie 
on any point interior to that ellipse. $\hfill\blacksquare$\\[-0.4cm]

\section{Generalizations, Extensions and Applications}
\label{GeneralizationsExtensionsandApplications}

By considering in (\ref{generalizationotthergegen})
the substitution $\nu=d/2-1$ 
and the map $\nu\mapsto -\nu/2$,
one obtains the formula
\begin{eqnarray*}
\frac
{\rho^{(d-1)/2}(1-\rho^2)^{\nu-(d-1)/2}}
{(1+\rho^2-2\rho x)^\nu}
&=&\frac{e^{-i\pi(\nu-(d-1)/2)}\Gamma(\frac{d-2}{2})}{2\sqrt{\pi}\,\Gamma(\nu)}
\nonumber\\[0.0cm]
&&\hspace{-0.2cm}
\times\sum_{{n}=0}^\infty(2{n}+d-2)
Q_{{n}+(d-3)/2}^{\nu-(d-1)/2}
\left(\frac{1+\rho^2}{2\rho}\right)
C_{n}^{d/2-1}(x).
\end{eqnarray*}
This formula
generalizes (9.9.2) in Andrews, Askey \& Roy (1999) \cite{AAR}.  

\medskip

By taking the limit as $\mu\to 1/2$ in (\ref{biggeneralizationgegen2}), 
one obtains a general result is an expansion over Legendre polynomials, namely
\begin{equation}
\frac{1}{(z-x)^\nu}=\frac{e^{i\pi(1-\nu)}(z^2-1)^{(1-\nu)/2}}{\Gamma(\nu)}
\sum_{n=0}^\infty(2n+1)Q_n^{\nu-1}(z)P_n(x),
\label{generlizationofheinesformula}
\end{equation}
using 
\begin{equation}
P_n(x)=C_n^{1/2}(x),
\label{PlCl}
\end{equation}
which is clear by comparing 
(cf.~(18.7.9) of Olver {\it et al.} (2010) and (\ref{gegpolydefn}) or 
(\ref{Jacobidefn}))
\begin{equation}
P_{n}(x):={}_2F_1\left(-{n},{n}+1;1;\frac{1-x}{2}\right),
\label{legendrepolydefn}
\end{equation}
and (\ref{gegpolydefn}).
If one takes $\nu=1$ in (\ref{generlizationofheinesformula})
then one has an expansion of the Cauchy denominator
which generalizes Heine's formula 
(see for instance Olver (1997) \cite[Ex.~13.1]{Olver};
Heine (1878) \cite[p.~78]{Heine1878})
\[
\frac{1}{z-x}=\sum_{n=0}^\infty(2n+1)Q_n(z)P_n(x).
\]

By taking the limit as $\mu\to 1$ in (\ref{biggeneralizationgegen2}), 
one obtains a general result which is an expansion over Chebyshev polynomials
of the second kind, namely
\begin{equation}
\frac{1}{(z-x)^\nu}=\frac{2^{3/2}e^{i\pi(3/2-\nu)}}
{\sqrt{\pi}\,\Gamma(\nu)(z^2-1)^{\nu/2-3/4}}
\sum_{n=0}^\infty(n+1)Q_{n+1/2}^{\nu-3/2}(z)U_n(x),
\label{chebyshevsecondform}
\end{equation}
using (18.7.4) in Olver {\it et al.} (2010) \cite{NIST},
$U_n(x)=C_n^{1}(x).$
If one considers the case $\nu=1$ in (\ref{chebyshevsecondform}) then 
the associated Legendre function of the second kind reduces to an elementary
function through (8.6.11) in Abramowitz \& Stegun (1972) \cite{Abra}, namely
\[
\frac{1}{z-x}=2\sum_{n=0}^\infty \frac{U_n(x)}{(z+\sqrt{z^2-1})^{n+1}}.
\]
By taking the limit as $\nu\to 1$ in (\ref{biggeneralizationgegen2}), 
one produces the Gegenbauer expansion of the Cauchy denominator given in 
Durand, Fishbane \& Simmons (1976) \cite[(7.2)]{DurandFishSim}, 
namely
\[
\frac{1}{z-x}=\frac{2^{\mu+1/2}}{\sqrt{\pi}}\Gamma(\mu)e^{i\pi(\mu-1/2)}(z^2-1)^{\mu/2-1/4}
\sum_{n=0}^\infty(n+\mu)Q_{n+\mu-1/2}^{-\mu+1/2}(z)C_n^\mu(x).
\]
Using (2.4) therein, the associated Legendre function of the second kind
is converted to the Gegenbauer function of the second kind. 

If one considers the limit as $\mu\to 0$ in (\ref{biggeneralizationgegen2})
using
\[
\lim_{\mu\to 0}\frac{{n}+\mu}{\mu}C_{n}^\mu(x)=\epsilon_{n} T_{n}(x)
\]
(see for instance (6.4.13) in Andrews, Askey \& Roy (1999) \cite{AAR}),
where $T_{n}:\C\to\C$ is the Chebyshev polynomial of the first kind 
defined as (see Section 5.7.2 in Magnus, Oberhettinger \& Soni (1966) \cite{MOS})
\[
T_n(x):=
{}_2F_1\left(-n,n;\frac12;\frac{1-x}{2}\right),
\]
and $\epsilon_{n}=2-\delta_{{n},0}$ is the Neumann factor, commonly 
appearing in Fourier cosine series, then one obtains
\begin{equation}
\frac{1}{(z-x)^\nu}=\sqrt{\frac{2}{\pi}}\frac{e^{-i\pi(\nu-1/2)}(z^2-1)^{-\nu/2+1/4}}
{\Gamma(\nu)}
\sum_{n=0}^\infty \epsilon_nT_n(x) Q_{n-1/2}^{\nu-1/2}(z).
\label{recipchebyshev}
\end{equation}
The result (\ref{recipchebyshev}), which is a generalization of Heine's reciprocal 
square-root identity 
(see Heine (1881) \cite[p.~286]{Heine}; Cohl \& Tohline (1999) \cite[(A5)]{CT}).
Polynomials in $(z-x)$ also naturally arise by considering the limit 
$\nu\to n\in-\N_0.$  This limit is given in (4.4) of Cohl \& Dominici (2010) \cite{CohlDominici},
namely
\begin{equation}
(z-x)^q=i(-1)^{q+1}\sqrt{\frac{2}{\pi}}(z^2-1)^{q/2+1/4}\sum_{n=0}^q \epsilon_n T_n(x)\frac{(-q)_n}{(q+n)!}Q_{n-1/2}^{q+1/2}(z),
\label{polyzmx}
\end{equation}
for $q\in\N_0$. 
\medskip

\noindent Note that all of the above formulae are restricted by the convergence criterion
given by {\sc Theorem} 9.1.1 in Szeg\H{o} (1959) \cite{Szego}
(Expansion of an analytic function in a Jacobi series), i.e.,
since the functions on the left-hand side are analytic in $[-1,1]$, then 
the expansion formulae are convergent if the point $z\in\C$ lies on any 
ellipse with foci at $\pm 1$ then $x$ can lie on any point interior 
to that ellipse.  Except of course (\ref{polyzmx}) which converges for
all points $z,x\in\C$ since the function on the left-hand side is entire.

\medskip
An interesting extension of the results presented in this paper, originally
uploaded to {\tt arXiv} in Cohl (2011) \cite{CohlGengegenarXiv} have
been obtained recently in Szmytkowski (2011) \cite{SzmyGeg}, 
to obtain formulas such as
\begin{eqnarray}
\hspace{-0.1cm}
\sum_{n=0}^\infty\frac{n+\mu}{\mu}\mathrm{P}_{n+\mu-1/2}^{\nu-\mu}(t)C_n^\mu(x)=
\frac{\sqrt{\pi}\,(1-t^2)^{(\nu-\mu)/2}}
{2^{\mu-1/2}\Gamma(\mu+1)\Gamma\left(\frac12-\nu\right)}
\nonumber\\[0.2cm]
&&\hspace{-4.0cm}\times
\left\{ \begin{array}{ll}
\displaystyle 0 & \qquad\mathrm{if}\ -1<x<t<1,\\[5pt]
\displaystyle (x-t)^{-\nu-1/2} & \qquad\mathrm{if}\ -1<t<x<1, \nonumber
\end{array} \right.
\end{eqnarray}
and
\begin{eqnarray}
\hspace{-0.5cm}
\sum_{n=0}^\infty\frac{n+\mu}{\mu}\mathrm{Q}_{n+\mu-1/2}^{\nu-\mu}(t)C_n^\mu(x)=
\frac{\sqrt{\pi}\,\Gamma\left(\nu+\frac12\right)(1-t^2)^{(\nu-\mu)/2}}
{2^{\mu+1/2}\Gamma(\mu+1)}
\nonumber\\[0.2cm]
&&\hspace{-7.0cm}\times
\left\{ \begin{array}{ll}
\displaystyle (t-x)^{-\nu-1/2} & \qquad\mathrm{if}\ -1<x<t<1,\\[5pt]
\displaystyle (x-t)^{-\nu-1/2}\cos[\pi(\nu+\tfrac{1}{2})] & \qquad\mathrm{if}\ -1<t<x<1, \nonumber
\end{array} \right.
\end{eqnarray}
where $\mbox{Re}\,\mu>-1/2$, $\mbox{Re}\,\nu<1/2$ and
$\mathrm{P}_\nu^\mu,\mathrm{Q}_\nu^\mu:(-1,1)\to\C$ are Ferrers functions
(associated Legendre functions on-the-cut) of the first and second kind.
The Ferrers functions of the first and second kind 
can be defined using
Olver {\it et al.} (2010) \cite[(14.3.11-12)]{NIST}.

\section{Expansion of a power-law fundamental solution 
of the polyharmonic equation}
\label{dge3}

A fundamental solution for the polyharmonic equation
on Euclidean space $\R^d$ is a function 
$\mcg_k^d:(\R^d\times\R^d)\setminus\{(\bfx,\bfx):\bfx\in\R^d\}\to\R$ 
which satisfies the equation
\begin{equation}
(-\Delta)^k{\mcg}_k^d({\bf x},{\bf x}^\prime)=\delta({\bf x}-{\bf x}^\prime),
\label{unnormalizedfundsolnpolydefn}
\end{equation}
where $\Delta:C^p(\R^d)\to C^{p-2}(\R^d),$ $p \ge 2$, is the Laplacian operator on $\R^d$ 
defined by
\[
\Delta:=\sum^d_{i=1}\frac{\partial^2}{\partial x_i^2},
\]
$\bfx=(x_1,\ldots,x_d), \bfxp=(x_1^\prime,\ldots,x_d^\prime)\in\R^d$, and
$\delta$ is the Dirac delta function.  Note that we introduce a minus sign into the 
equations where the Laplacian is used, such as in (\ref{unnormalizedfundsolnpolydefn}),
to make the resulting operator positive.
By Euclidean space $\R^d$, we mean the normed vector space given by the pair
$(\R^d,\|\cdot\|)$, where $\|\cdot\|:\R^d\to[0,\infty)$ is the Euclidean norm on $\R^d$ defined by
$\|\bfx\|:=\sqrt{x_1^2+\cdots+x_d^2},$
with inner product $(\cdot,\cdot):\R^d\times\R^d\to\R$ defined as
\begin{equation}
(\bfx,\bfxp):=\sum_{i=1}^d x_ix_i^\prime.
\label{eucinnerproduct}
\end{equation}
Then $\R^d$ is a $C^\infty$ Riemannian manifold with Riemannian metric 
induced from the inner product (\ref{eucinnerproduct}).  Set 
$\Si^{d-1}=\{\bfx\in\R^d:(\bfx,\bfx)=1\},$
then $\Si^{d-1},$ the $(d-1)$-dimensional unit hypersphere, 
is a regular submanifold of $\R^d$
and a $C^\infty$ Riemannian manifold with Riemannian metric induced from that
on $\R^d$.

\begin{theorem} Let $d,k\in\N$.  Define
\[
\mcg_k^d({\bf x},{\bf x}^\prime):=
\left\{ \begin{array}{ll}

{\displaystyle \frac{(-1)^{k+d/2+1}\ \|{\bf x}-{\bf x}^\prime\|^{2k-d}}
{(k-1)!\ \left(k-d/2\right)!\ 2^{2k-1}\pi^{d/2}}
\left(\log\|{\bf x}-{\bf x}^\prime\|-\beta_{p,d}\right)}\\[2pt]
\hspace{7.4cm} \mathrm{if}\  d\,\,\mathrm{even},\ k\ge d/2,\\[10pt]
{\displaystyle \frac{\Gamma(d/2-k)\|{\bf x}-{\bf x}^\prime\|^{2k-d}}
{(k-1)!\ 2^{2k}\pi^{d/2}}} \hspace{3.37cm} \mathrm{otherwise},
\end{array} \right.
\]
where $p=k-d/2$, $\beta_{p,d}\in\Q$ is defined as
$\beta_{p,d}:=\frac12\left[H_p+H_{d/2+p-1}-H_{d/2-1} \right],$
with $H_j\in\Q$ being the $j$th harmonic number
\[
H_j:=\sum_{i=1}^j\frac1i,
\]
then $\mcg_k^d$ is a fundamental solution for $(-\Delta)^k$
on Euclidean space $\R^d$.
\label{greenpoly}
\end{theorem}
\noindent {\it Proof }\  See Cohl (2010) \cite{CohlthesisII} and Boyling (1996) \cite{Boyl}.\\[0.2cm]
\medskip

Consider the following functions 
$\g_k^d,\li_k^d:(\R^d\times\R^d)\setminus\{(\bfx,\bfx):\bfx\in\R^d\}\to\R$ 
defined for $d$ odd and for $d$ even with $k\le d/2-1$ as a power-law, namely
\begin{equation}
\g_k^d(\bfx,\bfxp):=\|\bfx-\bfxp\|^{2k-d},
\label{frakgkd}
\end{equation}
and for $d$ even, $k\ge d/2$, with logarithmic behavior as
\[
\li_k^d(\bfx,\bfxp):=\|\bfx-\bfxp\|^{2p}
\left(\log\|\bfx-\bfxp\|-\beta_{p,d}\right),
\]
with $p=k-d/2$.  By {\sc Theorem} \ref{greenpoly} we see that 
the functions $\g_k^d$ and $\li_k^d$ equal real non-zero constant multiples
of $\mcg_k^d$ for appropriate parameters. Therefore by 
(\ref{unnormalizedfundsolnpolydefn}),
$\g_k^d$ and $\li_k^d$ are fundamental 
solutions
of the polyharmonic equation for appropriate parameters.  In this paper, 
we only consider functions with power-law behavior, although in future
publications we will consider the logarithmic case
(see Cohl (2012) \cite{Cohl12log} for the relevant Fourier expansions).

Now we consider the set of hyperspherical coordinate systems which parametrize points 
on $\R^d$.  The Euclidean distance between two points represented in these 
coordinate systems is given by
\[
\displaystyle \|\bfx-\bfxp\|=\sqrt{2 rr^\prime}
\left[z-\cos\gamma\right]^{1/2},
\]
where the toroidal parameter $z\in(1,\infty),$
(2.6) in Cohl {\it et al.} (2001) \cite{CRTB}, is given by
\[
z:=\frac{r^2+{r^\prime}^2}
{\displaystyle 2rr^\prime},
\]

\noindent and the separation angle $\gamma\in[0,\pi]$ is given through 
\begin{equation}
\cos\gamma=\frac{(\bfx,{\bf x^\prime})}{rr^\prime},
\label{separationanglegen}
\end{equation}

\noindent where $r,r^\prime\in[0,\infty)$ are defined such 
that $r:=\|{\bf x}\|$ and $r^\prime:=\|{\bf x^\prime}\|.$
We will use these quantities to derive Gegenbauer expansions of 
power-law fundamental solutions for powers of the Laplacian $\g_k^d$ 
(\ref{frakgkd}) represented in hyperspherical coordinates.

\begin{corollary}
For $d\in\{3,4,5,\ldots\}$, $\nu\in\C$, $\bfx,\bfxp\in\R^d$ with $r=\|\bfx\|$, $r^\prime=\|\bfxp\|$, 
and $\cos\gamma=(\bfx,\bfxp)/(rr^\prime)$, the following formula holds
\begin{eqnarray}
\hspace{-7cm}\|\bfx-\bfxp\|^\nu&=&
\frac{e^{i\pi(\nu+d-1)/2}
\Gamma\left(\frac{d-2}{2}\right)}
{2\sqrt{\pi}\,\Gamma\left(-\frac{\nu}{2}\right)}
\frac{\left(r_>^2-r_<^2\right)^{(\nu+d-1)/2}}
{\left(rr^\prime\right)^{(d-1)/2}}
\nonumber\\[0.0cm]
&&\hspace{-0.00cm}\times\sum_{{n}=0}^\infty
\left(2{n}+d-2\right)
Q_{{n}+(d-3)/2}^{(1-\nu-d)/2}
\biggl(\frac{r^2+{r^\prime}^2}{2rr^\prime}\biggr)
C_{n}^{d/2-1}(\cos\gamma),
\label{expandgegenpowq}
\end{eqnarray}
where $r_\lessgtr:={\min \atop \max}\{r,r^\prime\}$.
\label{MYCOROLLARY}
\end{corollary}
\noindent Note that (\ref{expandgegenpowq})
is seen to be a generalization of Laplace's expansion on $\R^3$ 
(see for instance Sack (1964) \cite{Sacka})
\[
\frac{1}{\|\bfx-\bfxp\|}=\sum_{n=0}^\infty \frac{r_<^n}{r_>^{n+1}}P_n(\cos\gamma),
\]
which is demonstrated by utilizing (\ref{PlCl}) and
simplifying the 
associated Legendre function of the second kind in (\ref{expandgegenpowq}) through 
$Q_{-1/2}^{1/2}:\C\setminus(-\infty,1]\to\C$ defined such that
\[
Q_{-1/2}^{1/2}(z)=i\sqrt{\frac{\pi}{2}}(z^2-1)^{-1/4}
\]
(cf.~(8.6.10) in Abramowitz \& Stegun (1972) \cite{Abra}
and (\ref{assoclegQseckinddefn2F1})).

The {\it addition theorem 
for hyperspherical harmonics}, which generalizes
\begin{equation}
P_n(\cos\gamma)=\frac{4\pi}{2{n}+1}\sum_{m=-{n}}^{n}
Y_{{n},m}({\mathbf {\widehat x}})
\overline{Y_{{n},m}
(\mathbf {\widehat x^\prime)}},
\label{addtheorem}
\end{equation}
where $P_{n}(x)$ is the Legendre polynomial
of degree ${n}\in\N_0$,
for $d=3$, is given by
\begin{equation}
\sum_{K}
Y_n^K
(\wbfx)
\overline{Y_n^K
(\wbfxp)
}
=\frac{\Gamma(d/2)}{2\pi^{d/2}(d-2)}
(2n+d-2)
C_n^{d/2-1}(\cos\gamma),
\label{additionthmhypspheharm}
\end{equation}
where $K$ stands for a set of $(d-2)$-quantum 
numbers identifying degenerate harmonics
for a given value of ${n}$ and $d$, and
$\gamma$ is the separation angle (\ref{separationanglegen}).
The functions
$Y_{{n}}^{K}:\Si^{d-1}\to\C$ are the normalized hyperspherical harmonics,
and $Y_{n,k}:\Si^2\to\C$ are the normalized spherical harmonics for $d=3$.
Note that 
${\mathbf {\widehat x}},{\mathbf {\widehat x^\prime}}\in\Si^{d-1}$,
are the vectors of unit length in the direction of $\bfx, \bfxp \in \R^d$ 
respectively.  
For a proof of the addition theorem for hyperspherical harmonics
(\ref{additionthmhypspheharm}), see Wen \& Avery (1985) \cite{WenAvery} and for a relevant discussion, 
see Section 10.2.1 in Fano \& Rau (1996) \cite{FanoRau}.  
The correspondence between (\ref{addtheorem}) and (\ref{additionthmhypspheharm}) arises
from (\ref{gegpolydefn}) and (\ref{legendrepolydefn}) namely (\ref{PlCl}).

One can use the addition theorem for hyperspherical harmonics
to expand a fundamental solution of the polyharmonic equation on $\R^d$.
Through the use of the addition theorem for hyperspherical harmonics we see 
that the Gegenbauer polynomials $C_{n}^{d/2-1}(\cos\gamma)$ are hyperspherical
harmonics when regarded as a function of $\wbfx$ only (see 
Vilenkin (1968) \cite{Vilen}).
Normalization of the hyperspherical harmonics is accomplished through the following integral
\[
\int_{\Si^{d-1}}
Y_{n}^K
({\mathbf {\widehat x}})
\overline{Y_{n}^K
({\mathbf {\widehat x}})}
d\Omega
=1,
\]
\noindent where $d\Omega$ is the Riemannian volume measure on $\Si^{d-1}$.
The degeneracy, i.e., number of linearly independent solutions for 
a particular value of 
${n}$ and $d$, for the space of hyperspherical harmonics is given by
\begin{equation}
(2{n}+d-2)\frac{(d-3+{n})!}{{n}!(d-2)!}
\label{degeneracynumber}
\end{equation}
(see (9.2.11) in Vilenkin (1968) \cite{Vilen}).  
The total number of 
linearly independent solutions 
(\ref{degeneracynumber}) can be determined by counting the total number
of terms in the sum over $K$ in 
(\ref{additionthmhypspheharm}).
Note that this formula (\ref{degeneracynumber}) reduces to the standard result in 
$d=3$ with a degeneracy given by $2{n}+1$ and in $d=4$ with a degeneracy 
given by $({n}+1)^2$.

One can show the consistency of 
{\sc Corollary} \ref{MYCOROLLARY}
with the result for $d=2$ 
given by
\[
\displaystyle \|\bfx-\bfxp\|^\nu=
\frac{e^{i\pi(\nu+1)/2}}{\Gamma(-\nu/2)}
\frac{(r_>^2-r_<^2)^{(\nu+1)/2}}{\sqrt{\pi rr^\prime}}
\sum_{m=-\infty}^\infty
e^{im(\phi-\phi^\prime)}
Q_{m-1/2}^{-(\nu+1)/2}
\left(\frac{r^2+{r^\prime}^2}{2rr^\prime}\right),
\]
where 
$\nu\in\C\setminus\{0,2,4,\ldots\}$,
by considering the limit as $\mu\to 0$
in (\ref{biggeneralizationgegen2})
(see (\ref{recipchebyshev}) above).
These expansions are useful in that they allow one to perform azimuthal 
Fourier and Gegenbauer polynomial analysis for 
power-law fundamental solutions of the polyharmonic equation on $\R^d$.

\section{Conclusion}

In this paper, we introduced a generalization of the generating function for Gegenbauer 
polynomials which allows one to expand arbitrary powers of the distance
between two points on $d$-dimensional Euclidean space $\R^d$ in terms of hyperspherical
harmonics.   
This result has already found
physical applications such as in Szmytkowski (2011) \cite{Szmy6} who uses this
result to obtain a solution of the momentum-space Schr\"{o}dinger equation
for bound states of the $d$-dimensional Coulomb problem.
The Gegenbauer expansions presented in this 
paper can be used in conjunction with corresponding Fourier expansions
(Cohl \& Dominici (2010) \cite{CohlDominici}) to generate infinite sequences
of addition theorems for the Fourier coefficients
(see Cohl (2010) \cite{CohlthesisII}) of these expansions.
In future publications, we will present some of these addition theorems 
as well as extensions related to Fourier and Gegenbauer expansions
for logarithmic fundamental solutions of the polyharmonic equation.

\section*{Acknowledgements}
Much thanks to T.~H.~Koornwinder for valuable discussions.
Part of this work was conducted while H.~S.~Cohl was a National Research Council
Research Postdoctoral Associate in the Applied and Computational
Mathematics Division at the
National Institute of Standards and Technology, Gaithersburg, Maryland, U.S.A.



\newpage


\begin{thebibliography}{10}

\bibitem{Abra}
M.~Abramowitz and I.~A. Stegun.
\newblock {\em Handbook of mathematical functions with formulas, graphs, and
  mathematical tables}, volume~55 of {\em National Bureau of Standards Applied
  Mathematics Series}.
\newblock U.S. Government Printing Office, Washington, D.C., 1972.

\bibitem{AAR}
G.~E. Andrews, R.~Askey, and R.~Roy.
\newblock {\em Special functions}, volume~71 of {\em Encyclopedia of
  Mathematics and its Applications}.
\newblock Cambridge University Press, Cambridge, 1999.

\bibitem{Boyl}
J.~B. Boyling.
\newblock Green's functions for polynomials in the {L}aplacian.
\newblock {\em Zeitschrift f\"ur Angewandte Mathematik und Physik},
  47(3):485--492, 1996.

\bibitem{CohlthesisII}
H.~S. {Cohl}.
\newblock {\em Fourier and {G}egenbauer expansions for fundamental solutions of
  the {L}aplacian and powers in {${\mathbf R}^d$} and {${\mathbf H}^d$}}.
\newblock PhD thesis, The University of Auckland, Auckland, New Zealand, 2010.
\newblock xiv+190 pages.

\bibitem{Cohl12log}
H.~S. {Cohl}.
\newblock {Fourier expansions for a logarithmic fundamental solution of the
  polyharmonic equation}.
\newblock {\em {\tt arXiv:1202.1811}}, 2012.

\bibitem{CohlGengegenarXiv}
H.~S. {Cohl}.
\newblock {On a generalization of the generating function for Gegenbauer
  polynomials}.
\newblock {\em {\tt arXiv:1105.2735}}, 2012.

\bibitem{CohlDominici}
H.~S. {Cohl} and D.~E. {Dominici}.
\newblock {Generalized Heine's identity for complex Fourier series of
  binomials}.
\newblock {\em Proceedings of the Royal Society A}, 467:333--345, 2011.

\bibitem{CRTB}
H.~S. Cohl, A.~R.~P. Rau, J.~E. Tohline, D.~A. Browne, J.~E. Cazes, and E.~I.
  Barnes.
\newblock Useful alternative to the multipole expansion of $1/r$ potentials.
\newblock {\em Physical Review A: Atomic and Molecular Physics and Dynamics},
  64(5):052509, Oct 2001.

\bibitem{CT}
H.~S. {Cohl} and J.~E. {Tohline}.
\newblock {A Compact Cylindrical Green's Function Expansion for the Solution of
  Potential Problems}.
\newblock {\em The Astrophysical Journal}, 527:86--101, December 1999.

\bibitem{DurandFishSim}
L.~Durand, P.~M. Fishbane, and L.~M. Simmons, Jr.
\newblock Expansion formulas and addition theorems for {G}egenbauer functions.
\newblock {\em Journal of Mathematical Physics}, 17(11):1933--1948, 1976.

\bibitem{ErdelyiHTFII}
A.~Erd{\'e}lyi, W.~Magnus, F.~Oberhettinger, and F.~G. Tricomi.
\newblock {\em Higher transcendental functions. {V}ol. {II}}.
\newblock Robert E. Krieger Publishing Co. Inc., Melbourne, Fla., 1981.

\bibitem{FanoRau}
U.~Fano and A.~R.~P. Rau.
\newblock {\em Symmetries in quantum physics}.
\newblock Academic Press Inc., San Diego, CA, 1996.

\bibitem{Gegenbauer1874}
L.~{Gegenbauer}.
\newblock {\"{U}ber einige bestimmte Integrale}.
\newblock {\em Sitzungsberichte der Kaiserlichen Akademie der Wissenschaften.
  Mathematische-Naturwissenschaftliche Classe.}, 70:433--443, 1874.

\bibitem{Gegenbauer1877}
L.~{Gegenbauer}.
\newblock {\"{U}ber die Functionen $C_n^\nu(x)$}.
\newblock {\em Sitzungsberichte der Kaiserlichen Akademie der Wissenschaften.
  Mathematische-Naturwissenschaftliche Classe.}, 75:891--896, 1877.

\bibitem{Gegenbauer1884}
L.~{Gegenbauer}.
\newblock {Zur theorie der functionen $C_n^\nu(x)$}.
\newblock {\em Denkschriften der Kaiserlichen Akademie der Wissenschaften zu
  Wien. Mathematisch-naturwissenschaftliche Classe}, 48:293--316, 1884.

\bibitem{Gegenbauer1888}
L.~{Gegenbauer}.
\newblock {\"{U}ber die Functionen $C_n^\nu(x)$}.
\newblock {\em Sitzungsberichte der Kaiserlichen Akademie der Wissenschaften.
  Mathematische-Naturwissenschaftliche Classe.}, 97:259--270, 1888.

\bibitem{Gegenbauer1893}
L.~{Gegenbauer}.
\newblock {Das Additionstheorem der Functionen $C_n^\nu(x)$}.
\newblock {\em Sitzungsberichte der Kaiserlichen Akademie der Wissenschaften.
  Mathematische-Naturwissenschaftliche Classe.}, 102:942--950, 1893.

\bibitem{Heine1878}
E.~Heine.
\newblock {\em Handbuch der {K}ugelfunctionen, {T}heorie und {A}nwendungen
  (volume 1)}.
\newblock Druck und Verlag von G. Reimer, Berlin, 1878.

\bibitem{Heine}
E.~Heine.
\newblock {\em Handbuch der {K}ugelfunctionen, {T}heorie und {A}nwendungen
  (volume 2)}.
\newblock Druck und Verlag von G. Reimer, Berlin, 1881.

\bibitem{Ismail}
M.~E.~H. Ismail.
\newblock {\em Classical and quantum orthogonal polynomials in one variable},
  volume~98 of {\em Encyclopedia of Mathematics and its Applications}.
\newblock Cambridge University Press, Cambridge, 2005.
\newblock With two chapters by Walter Van Assche, With a foreword by Richard A.
  Askey.

\bibitem{MOS}
W.~Magnus, F.~Oberhettinger, and R.~P. Soni.
\newblock {\em Formulas and theorems for the special functions of mathematical
  physics}.
\newblock Third enlarged edition. Die Grundlehren der mathematischen
  Wissenschaften, Band 52. Springer-Verlag New York, Inc., New York, 1966.

\bibitem{Olver}
F.~W.~J. Olver.
\newblock {\em Asymptotics and special functions}.
\newblock AKP Classics. A K Peters Ltd., Wellesley, MA, 1997.
\newblock Reprint of the 1974 original [Academic Press, New York].

\bibitem{NIST}
F.~W.~J. Olver, D.~W. Lozier, R.~F. Boisvert, and C.~W. Clark, editors.
\newblock {\em N{IST} handbook of mathematical functions}.
\newblock Cambridge University Press, Cambridge, 2010.

\bibitem{Sacka}
R.~A. Sack.
\newblock Generalization of {L}aplace's expansion to arbitrary powers and
  functions of the distance between two points.
\newblock {\em Journal of Mathematical Physics}, 5:245--251, 1964.

\bibitem{SriManocha}
H.~M. Srivastava and H.~L. Manocha.
\newblock {\em A treatise on generating functions}.
\newblock Ellis Horwood Series: Mathematics and its Applications. Ellis Horwood
  Ltd., Chichester, 1984.

\bibitem{Szego}
G.~Szeg{\H{o}}.
\newblock {\em Orthogonal polynomials}.
\newblock American Mathematical Society Colloquium Publications, Vol. 23.
  Revised ed. American Mathematical Society, Providence, R.I., 1959.

\bibitem{Szmy6}
R.~Szmytkowski.
\newblock Alternative approach to the solution of the momentum-space
  {S}chr\"odinger equation for bound states of the {$N$}-dimensional {C}oulomb
  problem.
\newblock {\em Annalen der Physik}, 524(6-7):345--352, 2012.

\bibitem{SzmyGeg}
R.~{Szmytkowski}.
\newblock {Some integrals and series involving the Gegenbauer polynomials and
  the Legendre functions on the cut $(-1,1)$}.
\newblock {\em Integral Transforms and Special Functions},
  {23}({11}):{847--852}, {2012}.

\bibitem{Vilen}
N.~Ja. Vilenkin.
\newblock {\em Special functions and the theory of group representations}.
\newblock Translated from the Russian by V. N. Singh. Translations of
  Mathematical Monographs, Vol. 22. American Mathematical Society, Providence,
  R. I., 1968.

\bibitem{WenAvery}
Z.~Y. Wen and J.~Avery.
\newblock Some properties of hyperspherical harmonics.
\newblock {\em Journal of Mathematical Physics}, 26(3):396--403, 1985.

\end{thebibliography}

\end{document}